\documentclass[12pt]{amsart}%
\usepackage{amsmath}
\usepackage{graphicx}%
\usepackage{amsfonts}%
\usepackage{amssymb}
\newtheorem{theorem}{Theorem}[section]
\theoremstyle{remark}
\newtheorem*{acknowledgements}{Acknowledgements}
\numberwithin{equation}{section}
\makeatletter
\renewcommand*\subjclass[2][2000]{%
  \def\@subjclass{#2}%
  \@ifundefined{subjclassname@#1}{%
    \ClassWarning{\@classname}{Unknown edition (#1) of Mathematics
      Subject Classification; using '2000'.}%
  }{%
    \@xp\let\@xp\subjclassname\csname subjclassname@#1\endcsname
  }%
}
\makeatother
\renewcommand{\subjclassname}{\textup{2000} Mathematics Subject Classification}
\input cyracc.def

\mathcode`\;="303B
\begin{document}
\title[Partial differential equations and Lie groups]{Some second-order partial differential equations associated with Lie groups}
\author{Palle E. T. Jorgensen}
\address{Department of Mathematics\\
The University of Iowa\\
14 MacLean Hall\\
Iowa City, IA 52242-1419\\
U.S.A.}
\email{jorgen@math.uiowa.edu}
\urladdr{http://www.math.uiowa.edu/\symbol{126}jorgen/}
\thanks{This research was partially supported by two grants from the U.S. National
Science Foundation, and by the Centre for Mathematics and its Applications
(CMA) at The Australian National University (ANU)}
\thanks{This paper is an expanded version of a lecture given by the author at the
National Research Symposium on Geometric Analysis and Applications at the ANU
in June of 2000.}
\dedicatory{To Derek Robinson on the occasion of his 65th birthday}

\begin{abstract}
In this note we survey results in recent research papers on the use of Lie
groups in the study of partial differential equations. The focus will be on
parabolic equations, and we will show how the problems at hand have solutions
that seem natural in the context of Lie groups. The research is joint with
D.W. Robinson, as well as other researchers who are listed in the references.

\end{abstract}
\maketitle

\section{\label{Int}Introduction}

When the Hamiltonian of a quantum-mechanical system is related to a Lie
algebra, it is often possible to use the representation structure of the Lie
algebra to decompose the Hilbert space of the quantum-mechanical system into
simpler (irreducible) pieces. For example, if a Hamiltonian commutes with the
generators of a Lie algebra, the Hilbert space of the system can be decomposed
into irreducibles of the Lie algebra, and the Lie algebra elements themselves
can be used as elements in a set of commuting observables.

We have aimed at making the present paper accessible to a wide audience of
non-specialists, stressing the general ideas and motivating examples, as
opposed to technical details.

The class of such Hamiltonians is quite large: see \cite{JoKl85} and
\cite{Jor88}. In this introduction we will review those Hamiltonians $H$ whose
interaction terms are polynomial in the position variables. Such Hamiltonians
are directly and naturally related to \emph{nilpotent} Lie algebras. The
nilpotent case is studied in Section \ref{Per}.

The spectrum of $H$ is obtained by decomposing the physical space on which the
Hamiltonian $H$ acts into irreducible representations of the underlying
nilpotent group. Sometimes this decomposition is decisive, as is the case with
a particle in a constant magnetic field, where the decomposition leads to a
harmonic-oscillator Hamiltonian. Sometimes the decomposition leads to a new
Hamiltonian that requires further analysis, as is the case with a particle in
a curved magnetic field.

The time evolution of the system is obtained by solving the heat equation of
the underlying nilpotent Lie group. By writing the Hamiltonian as a quadratic
sum of Lie-algebra elements and then using the representation of these
Lie-algebra elements arising from the regular representation, it is possible
to write $e^{-tH}$ as the convolution of a kernel (which is a solution of the
heat equation) with a representation acting on the physical Hilbert space; see
\cite{Jor88}.

The simplest case of this spectral picture is as follows: Consider a
nonrelativistic spinless particle of mass $m$ in an external magnetic field
$\mathbf{B}\left(  \mathbf{x}\right)  $. The Hamiltonian for such a system is
given by%
\begin{equation}
H=\frac{1}{2m}\left(  \mathbf{p}-\frac{e}{c}\mathbf{A}\right)  ^{2}%
,\label{eqJoKl85.2.1}%
\end{equation}
where $\mathbf{p}=\frac{h}{i}\nabla$ and $\mathbf{A}$ is the vector potential
satisfying $\mathbf{B}=\mathbf{\nabla}\times\mathbf{A}$. Consider the
commutators%
\begin{align}
\left[  p_{i}-\frac{e}{c}a_{i},p_{j}-\frac{e}{c}a_{j}\right]   &  =-\frac
{h}{i}\frac{e}{c}\varepsilon_{ijk}b_{k},\nonumber\\
\left[  p_{i}-\frac{e}{c}a_{i},b_{j}\right]   &  =\frac{h}{i}\frac{\partial
b_{j}}{\partial x_{i}}\equiv\frac{h}{i}b_{ij},\label{eqJoKl85.2.2}\\
\left[  p_{i}-\frac{e}{c}a_{i},b_{jk}\right]   &  =\frac{h}{i}\frac{\partial
b_{jk}}{\partial x_{i}}\equiv\frac{h}{i}b_{ijk},\nonumber\\
\vdots\qquad &  \qquad\qquad\vdots\qquad\qquad,\nonumber
\end{align}
where $\mathbf{A}=\left(  a_{1},a_{2},a_{3}\right)  $, $\mathbf{B}=\left(
b_{1},b_{2},b_{3}\right)  $, $\mathbf{x}=\left(  x_{1},x_{2},x_{3}\right)  $.
If $\mathbf{B}$ is a polynomial in $\mathbf{x}$, eventually the derivatives of
$\mathbf{B}$ will give zero, so that the set of commutators closes. The
resulting Lie algebra formed by real linear combinations of the elements
\begin{equation}
p_{i}-\frac{e}{c}a_{i},\;b_{i},\;b_{ij},\;\dots\label{eqIntNew.3}%
\end{equation}
is therefore a nilpotent Lie algebra, and the Hamiltonian (\ref{eqJoKl85.2.1})
is quadratic in the first three Lie algebra elements $X_{i}:=\left(
p_{i}-\frac{e}{c}a_{i}\right)  $, $i=1,2,3$, from the list (\ref{eqIntNew.3}).

We show further in \cite{JoKl85} and \cite{Jor88} that there is a unitary
representation $U$ of $G$ on $L^{2}\left(  \mathbb{R}^{3}\right)  $ such that%
\[
2mH=dU\,\left(  \sum_{i=1}^{3}\left(  p_{i}-\frac{e}{c}a_{i}\right)
^{2}\right)  .
\]
If there is a constant of motion for the Lie-algebra elements $p_{i}-\frac
{e}{c}a_{i}$, then $U$ is a direct integral over a corresponding spectral
parameter $\xi$. We then get $H=\int^{\oplus}d\xi\,H^{\left(  \xi\right)  }$
where $H$ has absolutely continuous spectrum, while each $H^{\left(
\xi\right)  }$ has purely discrete spectrum. If $\lambda_{0}\left(
\xi\right)  \leq\lambda_{1}\left(  \xi\right)  \leq\cdots$ is the spectrum of
$H^{\left(  \xi\right)  }$, then each $\xi\mapsto\lambda_{i}\left(
\xi\right)  $ is real analytic, and we get the following typical spectral picture.%

\[
\setlength{\unitlength}{144bp}%
\begin{picture}
(2.25,1.25)(-1.125,-0.125) \put(-1,-0.02){\includegraphics[bb=0 72 288
216,height=144bp,width=288bp] {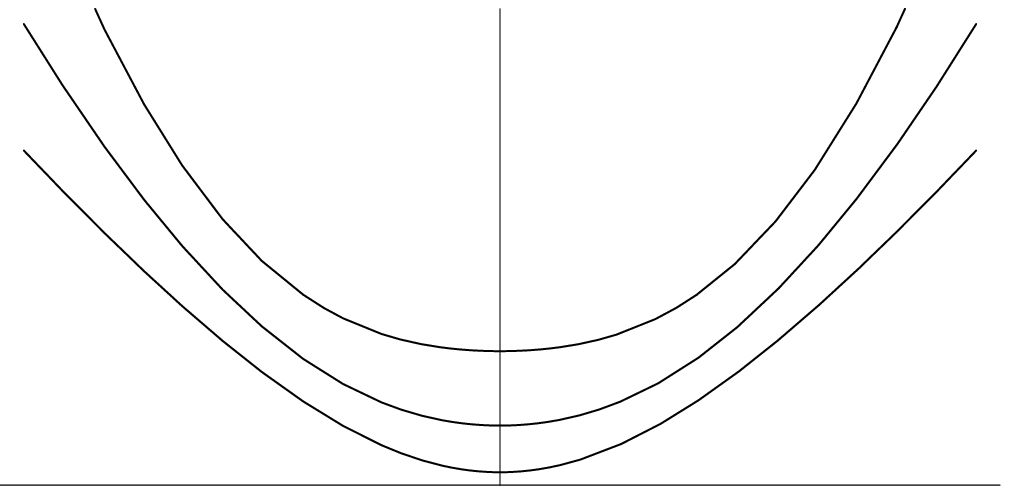}} \put(0.02,0.93){\makebox
(0,0)[l]{$\lambda$}} \put(0.95,-0.03){\makebox(0,0)[t]{$\xi$}} \put
(0.98,0.68){\makebox(0,0)[l]{$\lambda_{0}(\xi)$}} \put(0.98,0.945){\makebox
(0,0)[l]{$\lambda_{1}(\xi)$}} \put(0.8,1){\makebox(0,0)[b]{$\lambda_{2}(\xi)$%
}}
\end{picture}
\]

In this paper we will focus attention on a more restricted case wherein the
coefficients are periodic. As shown in Section \ref{GRd}, this case shares the
spectral band structure with the polynomial-magnetic-field case. We show that
in the periodic case the regularity of the coefficients may be relaxed, and in
fact, our spectral-theoretic results will be valid when the operator has
$L^{\infty}$-coefficients.

\section{\label{Per}Periodic operators}

We begin by recalling some elementary definitions and facts about stratified
Lie groups from \cite{FoSt82}. A real Lie algebra $\mathfrak{g}$ is
called\linebreak \ \emph{stratified} if it has a vector-space decomposition%
\begin{equation}
\mathfrak{g}=\bigoplus_{k=1}^{r}\mathfrak{g}^{\left(  k\right)  }%
,\label{eqPerNew.1}%
\end{equation}
for some $r$, which we shall take finite here, i.e., all but a finite number
of the subspaces $\mathfrak{g}^{\left(  k\right)  }$ are nonzero,%
\begin{equation}
\left[  \mathfrak{g}^{\left(  k\right)  },\mathfrak{g}^{\left(  l\right)
}\right]  \subseteq\mathfrak{g}^{\left(  k+l\right)  }\label{eqPerNew.2}%
\end{equation}
for all $k,l\in\mathbb{N}$, and $\mathfrak{g}^{\left(  1\right)  }$ generates
$\mathfrak{g}$ as a Lie algebra. Thus a stratified Lie algebra is
automatically nilpotent, and if $r$ is the largest integer such that
$\mathfrak{g}^{\left(  r\right)  }\neq0$, then $\mathfrak{g}$ is said to be
nilpotent of step $r$. A Lie group is defined to be \emph{stratified} if it is
connected and simply connected and its Lie algebra $\mathfrak{g}$ is stratified.

Let $G$ be a stratified Lie group and $\exp\colon\mathfrak{g}\rightarrow G$
the exponential map. The Campbell--Baker--Hausdorff formula establishes that%
\[
\exp\left(  X\right)  \exp\left(  Y\right)  =\exp\left(  H\left(  X,Y\right)
\right)  ,
\]
where $H\left(  X,Y\right)  =X+Y+\left[  X,Y\right]  /2+{}$a finite linear
combination of higher-order commutators in $X$ and $Y$. Thus $X,Y\rightarrow
H\left(  X,Y\right)  $ defines a group multiplication law on the underlying
vector space $V$ of $\mathfrak{g}$ which makes $V$ a Lie group whose Lie
algebra is $\mathfrak{g}$ and the exponential map $\exp\colon\mathfrak{g}%
\rightarrow V$ is simply the identity. Then $V$ with the group law is
diffeomorphic to $G$. Next let $d_{k}$ denote the dimension of $\mathfrak{g}%
^{\left(  k\right)  }$ and $d$ the dimension of $\mathfrak{g}$ and for each
$k$ choose a vector-space basis $X^{\left(  k\right)  }=\left(  X_{1}^{\left(
k\right)  },\dots,X_{d_{k}}^{\left(  k\right)  }\right)  $ of $\mathfrak{g}%
^{\left(  k\right)  }$ such that $X_{1},\dots,X_{d}=X_{1}^{\left(  1\right)
},\dots,X_{d_{r}}^{\left(  r\right)  }$ is a basis of $\mathfrak{g}$. If
$\xi_{1},\dots,\xi_{d}$ is the dual basis for $\mathfrak{g}^{\ast}$, i.e., if
$\xi_{k}\left(  X_{l}\right)  =\delta_{k,l}$, define $\eta_{k}=\xi_{k}%
\circ\exp^{-1}$. Then $\eta_{1},\dots,\eta_{d}$ are a system of global
coordinates for $G$, and the product rule on $G$ becomes%
\[
\eta_{k}\left(  xy\right)  =\eta_{k}\left(  x\right)  +\eta_{k}\left(
y\right)  +P_{k}\left(  x,y\right)  ,\qquad x,y\in G,
\]
where $P_{k}\left(  x,y\right)  $ is a finite sum of monomials in $\eta
_{i}\left(  x\right)  $, $\eta_{i}\left(  y\right)  $ for $i<k$ with degree
between $2$ and $m$. It follows that both left and right Haar measure on $G$
can be identified with Lebesgue measure $d\eta_{1}\,\cdots\,d\eta_{d}$.

If $X_{i}$ denotes one of the (abstract) Lie generators, we denote by $A_{i}$
the corresponding right-invariant vector field on $G$, i.e., $A_{i}$ on a test
function $\psi$ on $G$ is given by $A_{i}^{\left(  l\right)  }=dL\left(
X_{i}\right)  $, or more precisely,%
\begin{equation}
\left(  A_{i}^{\left(  l\right)  }\psi\right)  \left(  g\right)  =\frac
{d\;}{dt}\psi\left(  \exp\left(  -tX_{i}\right)  g\right)  |_{t=0},\qquad g\in
G,\label{eqPer.poundbis}%
\end{equation}
and similarly $A_{i}^{\left(  r\right)  }=dR\left(  X_{i}\right)  $ given by%
\begin{equation}
\left(  A_{i}^{\left(  r\right)  }\psi\right)  \left(  g\right)  =\frac
{d\;}{dt}\psi\left(  g\exp\left(  tX_{i}\right)  \right)  |_{t=0}%
.\label{eqPer.poundter}%
\end{equation}
Since we can pass from left to right with the adjoint representation, the
formulas may be written in one alone, and we will work with $A_{i}^{\left(
l\right)  }$, and denote it simply $A_{i}$.

If $1\leq j\leq d_{1}$ we will need the functions $y_{j}$ on $G$ defined by%
\begin{equation}
y_{j}\left(  \exp\left(  \sum_{k=1}^{d}\eta_{k}X_{k}\right)  \right)
=\eta_{j}. \label{eqPer.pound}%
\end{equation}
These functions satisfy the following system of differential equations:%
\begin{equation}
-A_{i}^{\left(  l\right)  }y_{j}=A_{i}^{\left(  r\right)  }y_{j}=\delta_{i,j}.
\label{eqBBJR95.3}%
\end{equation}
It follows by the standard ODE existence theorem that the functions $y_{i}$ on
$G$ are determined uniquely by (\ref{eqBBJR95.3}) and the ``initial''
conditions $y_{i}\left(  e\right)  =0$. Also note that (\ref{eqBBJR95.3}) is
consistent only for the differential equations defined from a sub-basis
$A_{1},\dots,A_{d_{1}}$, and that they would be overdetermined had we instead
used a basis: hence the distinction between subelliptic and elliptic.

In addition, we have given a discrete subgroup $\Gamma$ in $G$ such that
$M=G/\Gamma$ is compact. It is well-known that it then has a unique (up to
normalization) \cite{Jor88,Rob91} invariant measure $\mu$. The corresponding
Hilbert space is $L^{2}\left(  M,\mu\right)  $, and the invariant operators on
$G$ pass naturally to invariant operators on $M$; see \cite{BBJR95}. Let
$X_{1},\dots,X_{d_{1}}$ be the generating Lie-algebra elements. Then the
corresponding invariant vector fields on $G$ will be denoted $A_{1}%
,\dots,A_{d_{1}}$, and those on $M$ will be denoted $B_{1},\dots,B_{d_{1}}$.
Functions $c_{i,j}\in L^{\infty}\left(  G\right)  $ are given, and we form the
quadratic form%
\begin{equation}
h\left(  f\right)  =\sum_{i,j=1}^{d_{1}}\left\langle A_{i}f\mid c_{i,j}%
A_{j}f\right\rangle . \label{eqPer.1}%
\end{equation}
If further%
\begin{equation}
c_{i,j}\left(  g\gamma\right)  =c_{i,j}\left(  g\right)  \text{\qquad for
}g\in G,\;\gamma\in\Gamma, \label{eqPer.2}%
\end{equation}
then we have a corresponding form $h_{M}$ on $M=G/\Gamma$.

Introducing%
\begin{equation}
c_{i,j}^{\varepsilon}\left(  x\right)  =c_{i,j}\left(  \varepsilon
^{-1}x\right)  ,\qquad\varepsilon>0,\label{eqPer.3}%
\end{equation}
we get for each $\varepsilon$ a periodic problem corresponding to the period
lattice $\varepsilon\Gamma$. To speak about $\varepsilon\Gamma$ for
$\varepsilon\in\mathbb{R}_{+}$, we must have an action of $\mathbb{R}_{+}$ on
$G$ which generalizes the familiar one%
\[
\varepsilon\colon\left(  x_{1},\dots,x_{d}\right)  \longmapsto\left(
\varepsilon x_{1},\dots,\varepsilon x_{d}\right)
\]
of $\mathbb{R}^{d}$. It turns out that this can only be done if $G$ is
\emph{stratified}, and so in particular nilpotent; see \cite{FoSt82},
\cite{Jor88}. In that case it is possible to construct a group of
automorphisms $\left\{  \delta_{\varepsilon}\right\}  _{\varepsilon
\in\mathbb{R}_{+}}$ of $G$ which is determined by the differentiated action
$d\delta_{\varepsilon}$ on the Lie algebra $\mathfrak{g}$. If $\mathfrak{g}$
is specified as in (\ref{eqPerNew.1})--(\ref{eqPerNew.2}), then
\[
d\delta_{\varepsilon}\left(  X^{\left(  1\right)  }\right)  =\varepsilon
X^{\left(  1\right)  },\qquad X^{\left(  1\right)  }\in\mathfrak{g}^{\left(
1\right)  },\;\varepsilon\in\mathbb{R}_{+}.
\]

Let $H$, respectively $H_{\varepsilon}$, be the selfadjoint operators
associated to the period lattices $\Gamma$ and $\varepsilon\Gamma$ (see
\cite{BBJR95} or \cite{Rob91}), and let $S_{t}=e^{-tH}$, $S_{t}^{\varepsilon
}=e^{-tH_{\varepsilon}}$.

We now turn to the homogenization analysis of the limit $\varepsilon
\rightarrow0$ which leads to our comparison of the variable-coefficient case
to the constant-coefficient one. It should be stressed that in the Lie case,
even the ``constant-coefficient'' operator $\sum_{i,j}A_{i}\hat{c}_{i,j}A_{j}$
is not really constant-coefficient, as the vector fields $A_{i}$ are
variable-coefficient. 

Take even the simplest example where $G$ is the three-dimensional Heisenberg
group of upper triangular matrices of the form
\begin{equation}
g=%
\begin{pmatrix}
1 & x & z\\ 0 & 1 & y\\ 0 & 0 & 1
\end{pmatrix}
,\qquad x,y,z\in\mathbb{R}.\label{eqPer.dag}%
\end{equation}
In this case, $\dim\mathfrak{g}^{\left(  1\right)  }=2$, and $\dim$
$\mathfrak{g}^{\left(  2\right)  }=1$, with $\mathfrak{g}^{\left(  2\right)
}$ spanned by the central element in the Lie algebra. Differentiating matrix
multiplication (\ref{eqPer.dag}) on the left as in (\ref{eqPer.poundbis}), we
get the following three identities:%
\begin{align*}
A_{1}  & =\frac{\partial\;}{\partial x}+y\frac{\partial\;}{\partial
z}=-dL\left(  X_{1}\right)  ,\\
A_{2}  & =\frac{\partial\;}{\partial y}=-dL\left(  X_{2}\right)  ,\\
A_{3}  & =\frac{\partial\;}{\partial z}=-dL\left(  X_{3}\right)  ,
\end{align*}
where the first vector field is of course \emph{variable} coefficients.

We will use standard tools \cite{ZKO94} (see also \cite{Dau92}, \cite{Tho73}%
,\linebreak \ \cite{Wil78}) on homogenization.

\begin{theorem}
\label{ThmPer.1}\cite{BBJR95} Suppose the system $c_{i,j}\in L^{\infty}$ is
given and assumed strongly elliptic. Then there is a $C_{0}$-semigroup
$\hat{S}_{t}$ on $L^{2}\left(  G,dx\right)  $ with constant coefficients,
where $dx$ is left Haar measure, such that%
\begin{equation}
\lim_{\varepsilon\rightarrow0}\left\|  \left(  S_{t}^{\varepsilon}%
-\hat{S}_{t}^{{}}\right)  f\right\|  _{2}=0\label{eqPer.4}%
\end{equation}
for all $f\in L^{2}\left(  G,dx\right)  $ and $t>0$.
\end{theorem}

The constant coefficients of the limit operator $\hat{c}_{i,j}$ may be
determined as follows: We show in \cite{BBJR95} that if%
\begin{equation}
c_{i,j}\left(  g\right)  :=h\left(  g_{i}-y_{i},g_{j}-y_{j}\right)
\label{eqPer.5}%
\end{equation}
and if $C\left(  g\right)  $ is the corresponding quadratic form, then the
problem%
\begin{equation}
\inf_{g}C\left(  g\right)  =:\hat{C} \label{eqPer.6}%
\end{equation}
has a unique solution, i.e., the infimum is attained at $f_{1},\dots,f_{d_{1}%
}$ such that%
\begin{equation}
C\left(  f\right)  =\hat{C}. \label{eqPer.7}%
\end{equation}

The order relation which is used in the infimum consideration (\ref{eqPer.6})
is the usual order on hermitian matrices: For every $g$, the matrix $C\left(
g\right)  :=\left(  c_{i,j}\left(  g\right)  \right)  _{i,j=1}^{d_{1}}$ is
hermitian, and the matrix inequality $C\left(  g\right)  \geq\hat{C}$ may thus
be spelled out as follows:%
\[
\sum_{i,j}\bar{z}_{i}c_{i,j}\left(  g\right)  z_{j}\geq\sum_{i,j}\bar{z}%
_{i}\hat{c}_{i,j}z_{j}\text{\qquad for all }z_{1},\dots,z_{d_{1}}\in
\mathbb{C}.
\]
Solvability of this variational problem is part of the conclusion of our
analysis in \cite{BBJR95}, i.e., the existence of the minimizing functions
$f_{1},\dots,f_{d_{1}}$.

Then the coefficients of the homogenized operator can also be computed with
the aid of the coordinates $y_{i}$, $i=1,\dots,d_{1}$, introduced in
(\ref{eqPer.pound}) and (\ref{eqBBJR95.3}). One has the representation%
\begin{align}
\hat{c}_{i,j} &  =\int_{Y}dy\,\sum_{k,l=1}^{d_{1}}\left(  A_{k}\left(
f_{i}\left(  y\right)  -y_{i}\right)  \right)  c_{k,l}\left(  y\right)
\left(  A_{l}\left(  f_{j}\left(  y\right)  -y_{j}\right)  \right)
\label{eqBBJR95.11}\\
&  =h_{Y}\left(  f_{i}-y_{i},f_{j}-y_{j}\right)  ,\nonumber
\end{align}
where $h$ denotes the sesquilinear form associated with $H$, and the subscript
$Y$ refers to the region of integration. Specifically, $Y$ is a
\emph{fundamental domain} for the given lattice $\Gamma$ in $G$. For example,
we may take $Y$ to be defined by%
\begin{equation}
Y=\bigcap_{\gamma\in\Gamma}\left\{  x\in G;\left|  x\right|  \leq\left|
x\gamma^{-1}\right|  \right\}  ,\label{eqPer.gamma}%
\end{equation}
and $\left|  \,\cdot\,\right|  $ defined relative to a geodesic distance $d$,
$\left|  x\right|  :=d\left(  x,e\right)  $, $x\in G$. Then

\begin{enumerate}
\item \label{aftereqBBJR95.11(1)}$\bigcup_{\gamma\in\Gamma}Y\gamma=G$, and

\item \label{aftereqBBJR95.11(2)}$\operatorname*{meas}\left(  Y\gamma_{1}\cap
Y\gamma_{2}\right)  =0$ whenever $\gamma_{1}\neq\gamma_{2}$ in $\Gamma$.
\end{enumerate}

\noindent(These are the axioms for \emph{fundamental domains} of given
lattices, but we stress that (\ref{eqPer.gamma}) is just one choice in a vast
variety of possible choices.)

The simplest case of the construction is $G=\mathbb{R}$, and it was first
considered in \cite{Dav97} by Brian Davis. This is the simplest possible heat
equation, and we then have the conductivity represented by a periodic function
$c$, say%
\[
c\left(  x+p_{0}\right)  =c\left(  x\right)  ,\qquad x\in\mathbb{R},
\]
where $p_{0}$ is the period. Then $H=-\frac{d\;}{dx}c\left(  x\right)
\frac{d\;}{dx}$, and it can be checked that%
\[
\hat{c}=\left(  \frac{1}{p_{0}}\int_{0}^{p_{0}}\frac{dx}{c\left(  x\right)
}\right)  ^{-1}.
\]

\begin{theorem}
\label{ThmBBJR95.Pro4.2}\cite{BBJR95} Adopt the assumptions of Theorem
\textup{\ref{ThmPer.1}.} Then%
\[
\lim_{t\rightarrow\infty}t^{D/2}\operatorname*{ess\,sup}_{\left|  x\right|
^{2}+\left|  y\right|  ^{2}\leq at}\left|  K_{t}\left(  x;y\right)  -\hat
{K}_{t}\left(  x;y\right)  \right|  =0
\]
for each $a>0$ where $\left|  x\right|  =d_{c}\left(  x;e\right)  $, and
where
\begin{multline*}
d_{c}\left(  x;y\right)  =\sup\left\{  \psi\left(  x\right)  -\psi\left(
y\right)  ;\psi\in C_{c}^{\infty}\left(  G\right)  ,\vphantom{\sum
_i,j=1^d_1}\right. \\
\left.  \sum_{i,j=1}^{d_{1}}c_{i,j}\left(  A_{i}\psi\right)  \left(  A_{j}%
\psi\right)  \leq1\text{ pointwise}\right\}
\end{multline*}
and $A_{i}\psi$ refers to the Lie action of the vector field $A_{i}$ on $\psi$
from \textup{(\ref{eqPer.poundbis}).}
\end{theorem}

The number $D$ is the homogeneous degree defined from the given filtration
$\mathfrak{g}^{\left(  i\right)  }$ of the nilpotent Lie algebra
$\mathfrak{g}$. As spelled out in \cite{Jor88} and \cite{FoSt82}, there are
numbers $\nu_{i}$ depending on the Lie-structure coefficients such that%
\[
D=\sum_{i}\nu_{i}\dim\mathfrak{g}^{\left(  i\right)  }.
\]
To be specific, the numbers $\nu_{i}$ are determined in such a way that we get
a group of scaling automorphisms $\left\{  \delta_{\varepsilon}\right\}
_{\varepsilon\in\mathbb{R}_{+}}$ of $\mathfrak{g}$, and therefore on $G$, and
it is this group which is fundamental in the homogenization analysis.
Specifically, $\delta_{\varepsilon}\colon\mathfrak{g}\rightarrow\mathfrak{g}$
is defined by%
\begin{equation}
\delta_{\varepsilon}\left(  X^{\left(  i\right)  }\right)  =\nu_{i}X^{\left(
i\right)  },\qquad X^{\left(  i\right)  }\in\mathfrak{g}^{\left(  i\right)
},\label{eqPer.alpha}%
\end{equation}
and then extended to $\mathfrak{g}$ by linearity via (\ref{eqPerNew.1}), in
such a way that%
\begin{equation}
\delta_{\varepsilon}\left(  \left[  X,Y\right]  \right)  =\left[
\delta_{\varepsilon}\left(  X\right)  ,\delta_{\varepsilon}\left(  Y\right)
\right]  ,\qquad X,Y\in\mathfrak{g},\;\varepsilon\in\mathbb{R}_{+}%
.\label{eqPer.beta}%
\end{equation}
Hence if (\ref{eqPerNew.2}) holds, then it follows from (\ref{eqPer.alpha})
and (\ref{eqPer.beta}) that $\nu_{i}=i$ for $i=1,2,\dots$. In the case of the
Heisenberg Lie algebra $\mathfrak{g}$, we have $\left[  X,Y\right]  =Z$ as the
relation on the basis elements; $Z$ is central. Then $\mathfrak{g}^{\left(
1\right)  }=\operatorname*{span}\left(  X,Y\right)  $, $\mathfrak{g}^{\left(
2\right)  }=\mathbb{R}Z$, $\nu_{1}=1$, $\nu_{2}=2$, so $D=4$.

Let $K_{t}$ and $\hat{K}_{t}$ be the respective integral kernels for the
semigroups $S_{t}$ and $\hat{S}_{t}$, and set%
\[
\left|  \!\left|  \!\left|  K\right|  \!\right|  \!\right|  _{p}%
=\operatorname*{ess\,sup}_{x\in G}\left(  \int_{G}dy\,\left|  K\left(
x,t\right)  \right|  ^{p}\right)  ^{1/p}%
\]
and%
\[
\left|  \!\left|  \!\left|  K\right|  \!\right|  \!\right|  _{\infty
}=\operatorname*{ess\,sup}_{x,y\in G}\left|  K\left(  x,y\right)  \right|  .
\]
Then

\begin{theorem}
\label{ThmBBJR95Pro4.3}\cite{BBJR95} Adopt the assumptions of Theorem
\textup{\ref{ThmPer.1}.} Then%
\[
\lim_{t\rightarrow\infty}t^{D/2}\left|  \!\left|  \!\left|  K_{t}-\hat{K}%
_{t}\right|  \!\right|  \!\right|  _{\infty}=0,\qquad\lim_{t\rightarrow\infty
}\left|  \!\left|  \!\left|  K_{t}-\hat{K}_{t}\right|  \!\right|  \!\right|
_{1}=0.
\]
\end{theorem}

\section{\label{GRd}$G=\mathbb{R}^{d}$}

The case $G=\mathbb{R}^{d}$ was considered in \cite{BJR99}, where we further
showed that the limit $S_{t}^{\varepsilon}\rightarrow\hat{S}_{t}^{{}}$ then
holds also in the spectral sense. In that case, we scale by $\varepsilon=1/n$,
$n\rightarrow\infty$, and then identify the limit operator as having
absolutely continuous spectral type, and we prove spectral asymptotics. (A
general and classical reference for periodic operators is \cite{Eas73}.)

Starting with an equation which is invariant under the $\mathbb{Z}^{d}%
$-trans\-la\-tions, we then use the Zak transform \cite{Dau92} to write
$S_{t}=e^{-tH}$ as a direct integral over $\mathbb{T}^{d}$ ($=\mathbb{R}%
^{d}/\mathbb{Z}^{d}$), viz.,%
\begin{equation}
S_{t}=\int_{\mathbb{T}^{d}}^{\oplus}S_{t}^{\left(  z\right)  }%
,\label{eqGRd.poundpound}%
\end{equation}
and we establish continuity of $z\mapsto S_{t}^{\left(  z\right)  }$ in the
strong topology \cite[Lemma 2.2]{BJR99}. Pick a positive $C^{\infty}$-function
$\tau$ on $\mathbb{R}^{d}$ of integral one, and set
\[
c_{i,j}^{\left(  n\right)  }\left(  x\right)  =n^{d}\int_{\mathbb{R}^{d}%
}dy\,\tau\left(  ny\right)  c_{i,j}\left(  x-y\right)  ,
\]
and form the corresponding $C_{0}$-semigroup%
\[
S_{t}^{\left(  n\right)  }=e^{-tH^{\left(  n\right)  }},
\]
where $H^{\left(  n\right)  }$ is defined from $c_{i,j}^{\left(  n\right)  }$.
We then show in \cite{BJR99} that $S_{t}^{\left(  n\right)  }$ approximates
$S_{t}$, not only in the strong topology, but also in a spectral-theoretic
sense. Using this, we establish the following connection between
$S_{t}=e^{-tH}$ and $S_{t}^{\left(  z\right)  }=e^{-tH^{\left(  z\right)  }}$
in (\ref{eqGRd.poundpound}). Setting $z=\left(  e^{i\theta_{1}},\dots
,e^{i\theta_{d}}\right)  $, we get

\begin{theorem}
\label{ThmBJR99Cor4.7}If $\lambda_{n}\left(  z\right)  $ denotes the
eigenvalues of $H_{z}$ then%
\begin{multline}
\lim_{N\rightarrow\infty}\left\{  N^{2}\lambda_{n}\left(  w\right)
;w^{N}=z,\;n=0,1,\dots\right\}  \label{eqGRdNew.2}\\
=\left\{  \left\langle \left(  n-\theta\right)  \mid\hat{C}\left(
n-\theta\right)  \right\rangle ;n\in\mathbb{Z}^{d}\right\}  ,
\end{multline}
where the limit is in the sense of pointwise convergence of the ordered sets,
and where $\hat{C}=\left(  \hat{c}_{i,j}\right)  $ is the constant-coefficient
homogenized case.
\end{theorem}

The rate of convergence of the eigenvalues in (\ref{eqGRdNew.2}) can be
estimated further by a trace norm estimate.

We refer the reader to \cite{BJR99} for details of proof, but the arguments in
\cite{BJR99} are based in part on the references \cite{Aus96}, \cite{DaTr82},
\cite{Eas73}, and \cite{ZKON79}. In addition, we mention the papers
\cite{Aus96}, \cite{AMT98}, and \cite{TERo99}, which contain results which are
related, but with a different focus.

Finally, we mention that our result from \cite{BJR99}, Theorem
\ref{ThmBJR99Cor4.7}, has since been extended in several other directions:
see, e.g., \cite{Sob99} and \cite{She00}.

\begin{acknowledgements}
We are grateful to Brian Treadway for excellent typesetting and graphics
production, and to the participants in the National Research Symposium at The
Australian National University for fruitful discussions, especially A.F.M. ter Elst.
\end{acknowledgements}

\bibliographystyle{bftalpha}
\bibliography{jorgen}
\end{document}